\documentclass{article}
\usepackage{epsfig}
\usepackage{hyperref}
\usepackage{graphicx}
\usepackage{epstopdf}
\usepackage{amsmath}
\usepackage{setspace}
\usepackage{mathrsfs}
\usepackage{subfig}
\usepackage{color}
\usepackage{mathtools}
\usepackage{lipsum}
\usepackage[affil-it]{authblk}
\DeclarePairedDelimiter{\ceil}{\lceil}{\rceil}
\begin{document}

\title{Model-order reduction of lumped parameter systems via fractional calculus}

\author[1]{John P. Hollkamp}

\affil[1]{Department of Mechanical Engineering, Ray W. Herrick Laboratories, Purdue University, West Lafayette,
IN 47907
}

\author[2]{Mihir Sen}
\affil[2]{Department of Aerospace and Mechanical Engineering, University of Notre Dame, Notre Dame,
IN 46556
}
\author[1]{Fabio Semperlotti}

\date{Dated: \today}

\maketitle

\begin{abstract}
This study investigates the use of fractional order differential models to simulate the dynamic response of non-homogeneous discrete systems and to achieve efficient and accurate model order reduction. The traditional integer order approach to the simulation of non-homogeneous systems dictates the use of numerical solutions and often imposes stringent compromises between accuracy and computational performance. Fractional calculus provides an alternative approach where complex dynamical systems can be modeled with compact fractional equations that not only can still guarantee analytical solutions, but can also enable high levels of order reduction without compromising on accuracy. Different approaches are explored in order to transform the integer order model into a reduced order fractional model able to match the dynamic response of the initial system. Analytical and numerical results show that, under certain conditions, an exact match is possible and the resulting fractional differential models have both a complex and frequency-dependent order of the differential operator. The implications of this type of approach for both model order reduction and model synthesis are discussed.
\end{abstract}

\section{Introduction}
Numerical accuracy and computational efficiency have been long-standing challenges in the simulation of dynamical systems. Numerical solutions of complex continuous systems with non-trivial boundary and loading conditions typically require a discretization process to obtain lumped parameter models. The more complex the property spatial distribution (e.g. external loads, material or geometric parameters, boundary conditions), the higher the level of discretization needed to achieve a satisfactory representation of the original continuum system. The increased discretization directly impacts the computational performance and, for large systems, limits the level of achievable accuracy. This issue is even more accentuated when dealing with active control or real-time prediction of the dynamic response of systems under operating conditions for which fast state estimation is a key requirement. 

Over the last several decades, these challenges have motivated the rapid growth and development of methodologies for the synthesis of computationally efficient approaches able to reduce the overall size of the models (i.e. of the total number of Degrees Of Freedom - DOF) while maintaining high numerical accuracy and fidelity to the actual dynamics. These techniques, referred to as \textit{model order reduction}, are typically pursued when the dynamic response is sought only at selected locations (the so-called active DOFs) such that the DOFs associated with the remaining locations can be omitted. The reduction procedure is not trivial because it must account for the coupling between active and omitted DOFs in order to not change the underlying dynamics of the system. 

Many sources in the literature \cite{Besselink,deKlerk,Schilders} provide an extensive review of reduction techniques across a variety of disciplines. In the following, we concentrate only on applications to structural dynamics since this is the emphasis of our study. One widely used technique is Guyan reduction. This approach is also known as static reduction because it does not account for the system inertia and is therefore limited to statics. Reduction techniques for dynamic systems (therefore accounting for the system's inertia) are often based on mode superposition or component mode synthesis \cite{Besselink}. Perhaps one of the most widely used component mode synthesis technique, also known as the Craig-Bampton method \cite{CraigJA,CraigB}. The Craig-Bampton method divides the system into several substructures which are required to be compatible along their shared boundaries. Assuming these boundaries are held fixed, the Craig-Bampton method is able to combine the motion of these boundary points with the displacement modes of the substructures (known as constraint modes). The dynamics of the system can be reduced to a set of both fixed-interface and constraint modes \cite{Kuether}.

Many existing order reduction techniques can only provide an approximation of the local response. For instance, the accuracy of the approximation of the Craig-Bampton method is strictly dependent on the number of modes retained in the modal basis. The truncation of the basis should be assessed with respect to either the modal densities associated with the omitted modes or the dynamic content that should be transferred to the active degrees. In a similar way, enlarging the modal basis (i.e. extending the truncation order) comes at the expenses of computational performance. 

To address these limitations, we explored a reduction order technique based on fractional calculus. We will show that while fractional models contain less DOFs than the original system, often times the dynamic response can be matched exactly at the active degrees. This is an important advantage of our fractional order reduction over reduction order techniques typically used. In addition to the order reduction capabilities, we anticipate that the proposed approach has possible applications to the system identification based only on measured or experimental data. As an example of these capabilities, we will show the application of the fractional models to perform broadband system identification of discrete parameter systems.

While the mathematics of fractional calculus has been extensively studied in the past century, applications are relatively recent. In particular, fractional calculus has seen applications in engineering areas such as controls, visco- and thermo-elasticity \cite{Torvik,Wharmby,Achar1}, and wave propagation in complex media \cite{Fellah,Casasanta,Tarasov}. The reader is referred to \cite{Podlubny,Herrmann,Diethelm} for detailed reviews on the fundamentals of fractional calculus. In \S \ref{App_A} we define some basic fractional calculus quantities. We anticipate that the methodology discussed below produces fractional differential models of complex and frequency-dependent order. Love \cite{Love}, Ortigueira \cite{Ortigueira}, Ross \cite{Ross}, and Andriambololona \cite{Andriambololona} among others have developed the mathematics of complex fractional derivatives as well as potential uses. Adams \cite{Adams} and Neamaty \cite{Neamaty} have developed solutions to certain complex fractional order differential equations. Other authors including Atanackovic \cite{Atanackovi}, Makris \cite{Makris}, and Park \cite{Park} have successfully applied complex fractional calculus to viscoelasticity. While the mathematics of complex fractional calculus has been explored and developed, its connection to the actual physical processes being represented can be difficult to grasp. Perhaps Makris has one of the clearest interpretations of a complex order derivative: "\textit{...one may interpret the complex derivative of an arbitrary function as the superposition of complex derivatives of harmonic functions. Evidently, complex-order derivatives modulate the phase and amplitude of harmonic components of a time-dependent function in a more complicated way than real-order derivatives. An important difference between real-valued and complex-valued time derivatives is that phase modulation in the latter case is frequency dependent whereas in the former is not}" \cite{Makris}. As pointed out here by Makris, complex fractional derivatives can produce functions whose amplitude and phase are both frequency-dependent; this attribute plays a key role in the development of the fractional models. Valerio \cite{Valerio,ValerioJA} gives a thorough review of complex and variable-order fractional derivatives through the use of Laplace transforms and transfer functions.   

The remainder of the paper is structured as follows: in the next section, we present how to obtain an undamped fractional single degree of freedom model having the same dynamic response of an integer damped single degree of freedom. While this does not technically qualify as order reduction, it illustrates the basic methodology that will be used throughout this work. Then, we present the procedure to reduce a multiple integer-order degree of freedom system to a single fractional-order degree of freedom system. This represents the fundamental step of the order reduction methodology. Next, we extend the methodology to reduce a M-degree of freedom integer-order system to a N-degree of freedom fractional-order model, where $N < M$ and $M, N $. Lastly, we briefly discuss how the same methodology can be used in the frame of system identification and we show how to synthesize accurate fractional dynamic models based only on the knowledge of numerically obtained or experimentally measured data. 

An important remark should be made concerning the terminology used in this paper. It is well-known that, for integer order systems, the overall order of the system is strictly related to the number of degrees of freedom because each individual degree is assumed to behave as a second order system. Therefore, order and dimension are typically considered as equivalent concepts and used interchangeably. On the contrary, due to the infinite dimensional character of a fractional derivative, the connection between the overall order of a system and the number of its physical degrees of freedom is somewhat more ambiguous. Such a discussion goes well beyond the scope of this paper and we simply highlight that, in the remainder, the term \textit{degrees of freedom} will refer to the number of discrete masses while the term \textit{order} will refer to the order of the individual differential equations.

\section{Reduction to SDOF fractional systems}
This section discusses the fundamental approach to obtain a fractional undamped single degree of freedom (SDOF) model exhibiting a dynamic behavior equivalent to an integer damped single degree of freedom system. To facilitate the understanding of the procedure, we first discuss the conversion from an integer-order damped SDOF (I-SDOF) to a fractional-order SDOF (F-SDOF). Then, we will extend the methodology to the case of an integer-order multiple degree of freedom (I-MDOF) to be reduced to a F-SDOF.

\subsection{The fractional order SDOF model}

Consider the fractional SDOF oscillator in Figure \ref{D1}b. The equation of motion of the fractional model is given by:

\begin{figure}[ht] 
  \begin{center}
    \centerline{\includegraphics[scale=0.65,angle=0]{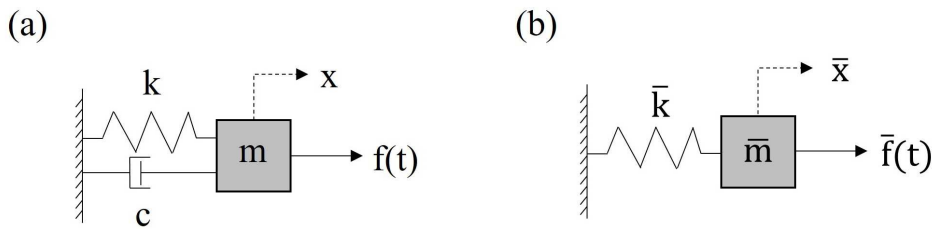}}
    \caption{Schematic of the (a) I-SDOF, and (b) the F-SDOF.}
    \label{D1}
  \end{center}
\end{figure}

\begin{eqnarray} \label{eq:FRO}
\bar{m}\frac{d^{\alpha(\omega)} \bar{x}}{dt^{\alpha(\omega)}}+\bar{k}\bar{x} = \bar{f}(t),
\end{eqnarray}
   
\noindent where $\bar{x}$ is the displacement from the equilibrium condition, $\bar{m}$ is the mass of the fractional model, $\bar{k}$ is the stiffness, $\bar{f}(t)$ is the time dependent load acting on the mass, and $\alpha(\omega)$ is the frequency-dependent order of the fractional derivative. In Equation \eqref{eq:FRO}, the fractional derivative should be intended in the Caputo form since this fractional derivative best represents an initial value problem of a physical system. However, it should be noted that the Riemann-Liouville derivative is equivalent to the Caputo derivative if all initial values are zero (which is also the general case considered in this paper). When 0 $< \alpha \leq$ 1, Equation \eqref{eq:FRO} is known as the fractional relaxation equation \cite{Podlubny,Herrmann}, while when 1 $< \alpha \leq$ 2, Equation \eqref{eq:FRO} is known as the fractional oscillation equation \cite{Podlubny,Herrmann}. Our model will use the form of Equation \eqref{eq:FRO}, but will not necessarily be restricted to 0 $< \alpha \leq$ 2. Also, our model will use a complex value of the fractional order $\alpha$ = $a+ib$, where $a$ and $b$ are the real and imaginary part of $\alpha$, respectively.

Taking the Laplace transform of Equation \eqref{eq:FRO} with zero initial conditions, we obtain:

\begin{eqnarray} \label{eq:LT1}
(\bar{m}s^{\alpha}+\bar{k})\bar{X}(s) = \bar{F}(s), 
\end{eqnarray}

\noindent where $s$ is the Laplace variable. The transfer function of the F-SDOF is:

\begin{eqnarray} \label{eq:TF}
G(s) = \frac{\bar{X}(s)}{\bar{F}(s)} = \bigg(\frac{1}{\bar{m}}\bigg)\frac{1}{s^\alpha + \frac{\bar{k}}{\bar{m}}}.
\end{eqnarray}

\noindent Substituting $\alpha = a+ib$ and $s=i\omega$, Equation \eqref{eq:TF} becomes:

\begin{eqnarray} \label{eq:TF2}
G(i\omega) = \bigg(\frac{1}{\bar{m}}\bigg)\frac{1}{(i\omega)^{a+ib}+\frac{\bar{k}}{\bar{m}}}.
\end{eqnarray}

\noindent Through some algebraic manipulation, Equation \eqref{eq:TF2} can be shown to be equivalent to: 

\begin{eqnarray} \label{eq:TF3}
G(i\omega) = \bigg(\frac{1}{\bar{m}}\bigg)\frac{1}{\gamma + \frac{\bar{k}}{\bar{m}}},
\end{eqnarray}

\noindent where: 

\begin{eqnarray} \label{eq:gamma}
\gamma = \omega^a e^{-\frac{b\pi}{2}} \omega^{ib} e^{i\frac{a\pi}{2}}.
\end{eqnarray}

\noindent In terms of sines and cosines, Equation \eqref{eq:gamma} can be written as:

\begin{equation} \label{eq:gamma2}
\begin{split}
\gamma = \omega^a e^{-\frac{b\pi}{2}}\bigg[& \textrm{cos}(\frac{a\pi}{2})\textrm{cos}(b \textrm{ln}(\omega))  -\textrm{sin}(\frac{a\pi}{2})\textrm{sin}(b \textrm{ln}(\omega))+ \\ 
&i\textrm{cos}(\frac{a\pi}{2})\textrm{sin}(b \textrm{ln}(\omega))+i\textrm{sin}(\frac{a\pi}{2})\textrm{cos}(b \textrm{ln}(\omega))\bigg].
\end{split}
\end{equation}

\noindent Equation \eqref{eq:TF3} can now be rewritten in terms of its real and imaginary parts:

\begin{eqnarray} \label{eq:TF4}
G(i\omega) = \bigg(\frac{1}{\bar{m}}\bigg)\frac{1}{\tau + i\xi},
\end{eqnarray}

\noindent where:

\begin{eqnarray} \label{eq:tau}
\tau = \frac{\bar{k}}{\bar{m}} + \omega^a e^{-\frac{b\pi}{2}}\bigg[\textrm{cos}(\frac{a\pi}{2})\textrm{cos}(b \textrm{ln}(\omega))-\textrm{sin}(\frac{a\pi}{2})\textrm{sin}(b \textrm{ln}(\omega))\bigg],
\end{eqnarray}

\begin{eqnarray} \label{eq:xi}
\xi = \omega^a e^{-\frac{b\pi}{2}}\bigg[\textrm{cos}(\frac{a\pi}{2})\textrm{sin}(b \textrm{ln}(\omega))+\textrm{sin}(\frac{a\pi}{2})\textrm{cos}(b \textrm{ln}(\omega))\bigg].
\end{eqnarray}

\noindent As a result, the magnitude and phase of the transfer function are: 

\begin{eqnarray} \label{eq:magn}
|G(i \omega)| = M = \frac{1}{\bar{m}\sqrt{\tau^2 + \xi^2}},
\end{eqnarray}

\begin{eqnarray} \label{eq:phase}
\angle G(i \omega) = \psi = -\textrm{tan}^{-1} \bigg(\frac{\xi}{\tau}\bigg).
\end{eqnarray}

If the forcing function is given by a harmonic load $F_0 \textrm{sin}(\omega t)$, then the steady state response of the F-SDOF is given by: 

\begin{eqnarray} \label{eq:FSS}
x(t) = F_0 M \textrm{sin}(\omega t + \psi).
\end{eqnarray}

\subsection{Exact conversion from a damped I-SDOF to a F-SDOF}\label{ISDOFtoFSDOF}
The general process to obtain the fractional derivative value $\alpha$ from a non-homogeneous integer order system is first illustrated by converting an I-SDOF into a F-SDOF having an equivalent dynamic response. 

Consider the differential equation of a classical damped oscillator (Figure \ref{D1}a):

\begin{eqnarray} \label{eq:ISDOF}
m\ddot{x}+c\dot{x}+kx = f(t),
\end{eqnarray}

\noindent where $m$ is the mass, $c$ is the damping coefficient, $k$ is the spring stiffness, $x$ is the displacement of the mass from equilibrium, and $f$ is the external force.

Assuming the forcing to be the same in the fractional model as in the integer model (that is, $f(t)$ = $\bar{f}(t)$), we can set the Laplace transform of Equation \eqref{eq:ISDOF} to be equal to the Laplace transform of the fractional model. Letting the displacements of the I-SDOF and the F-SDOF be equivalent (that is, $x$ = $\bar{x}$), we obtain a polynomial equation of fractional order:

\begin{eqnarray} \label{eq:LP}
ms^2+cs+k = \bar{m}s^\alpha + \bar{k}.
\end{eqnarray}

\noindent Since we are converting a SDOF to another SDOF, we set the coefficients of the fractional model equal to their integer model counterparts; that is, $\bar{m}$ = $m$ and $\bar{k}$ = $k$. Substituting this into Equation \eqref{eq:LP} and solving for $\alpha$ results in: 

\begin{eqnarray} \label{eq:alpha}
\alpha = 1 + \textrm{log}_{s}(s + \frac{c}{m}) = 1 + \frac{\textrm{ln}(s + \frac{c}{m})}{\textrm{ln}(s)}.
\end{eqnarray} 

\noindent Substituting s = i$\omega$ into Equation \eqref{eq:alpha} yields:

\begin{eqnarray} \label{eq:alpha2}
\alpha = 1 + \frac{\textrm{Ln}(i\omega + \frac{c}{m})}{\textrm{Ln}(i\omega)},
\end{eqnarray} 

\noindent where $\omega$ is the forcing frequency and Ln is the complex logarithm function. Recall that the complex logarithm function Ln$(z)$ in Equation \eqref{eq:alpha2} is defined as:

\begin{eqnarray}
\textrm{Ln}(z) = \textrm{ln}|z| + i\textrm{Arg}(z),
\end{eqnarray}

\noindent where $|z|$ is the modulus ($|z| = \sqrt{(\textrm{Re}(z))^2 + (\textrm{Im}(z))^2}$) and Arg$(z)$ is the principle argument (Arg$(z) = \textrm{arctan}(\textrm{Im}(z)/\textrm{Re}(z))$).

Equation \eqref{eq:alpha2} allows an important observation: the order $\alpha$ of the fractional model is both a complex and frequency-dependent quantity. Note that both characteristics are due to the fact that the damping term of the integer order model is now included into the fractional derivative; this is a well-known property of fractional oscillators \cite{Herrmann,Achar2,Ryabov,Stanislavsky,Tofighi}. As expected, setting $c=0$ would return an integer order derivative.  

\begin{figure}[!ht] 
  \begin{center}
    \centerline{\includegraphics[scale=0.7,angle=0]{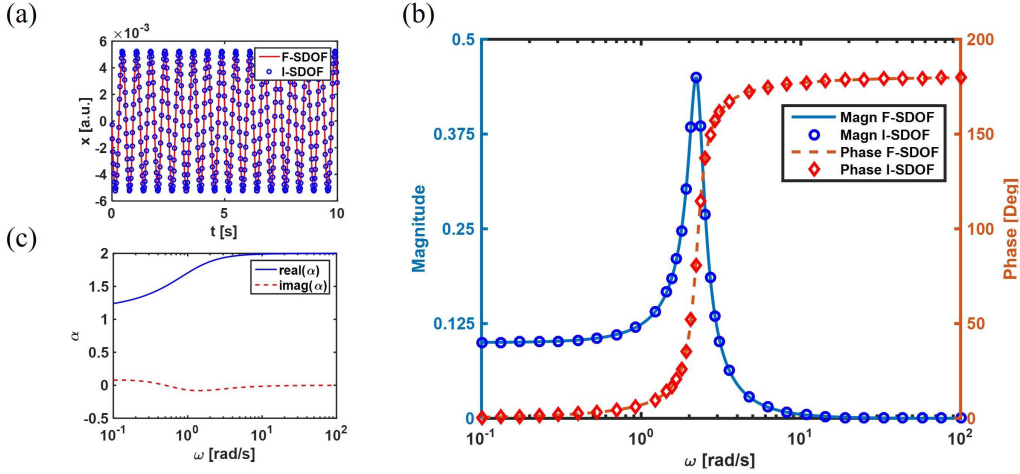}}
    \caption {a) Numerical simulations of the forced response showing the correspondence between the I-SDOF and the F-SDOF after the matching procedure is applied. Displacement time history at steady state for both the I-SDOF and the F-SDOF at a forcing frequency of 10 rad/s. b) Magnitude and phase of the transfer functions of I-SDOF and F-SDOF. c) Plot of the fractional order $\alpha(\omega)$ for the F-SDOF. }
    \label{SS}
  \end{center}
\end{figure}

The above results suggest that if $\alpha$ in Equation \eqref{eq:FRO} is chosen according to Equation \eqref{eq:alpha2}, then the steady state response of the two systems are exactly equivalent. To confirm this, we calculated numerically the response of the two systems. The steady state response of the I-SDOF subject to a force $F_o \textrm{sin}(\omega t)$ is given by $x(t) = Xsin(\omega t - \phi)$ \cite{Rao} while the steady state of the fractional model is given by Equation \eqref{eq:FSS}. Consider an I-SDOF where $m$ = 2, $c$ = 1, and $k$ = 10 (unit system is arbitrary, although we shall use units of seconds for time so to conveniently refer to frequency in terms of rad/s or Hz). As already discussed, we then let $\bar{m}$ = $m$ = 2 and $\bar{k}$ = $k$ = 10 in the fractional model. The value of $\alpha$ can then be calculated using Equation \ref{eq:alpha2} for a specified forcing frequency. The steady state response of the fractional model can be obtained using the magnitude and phase as given by Equations \ref{eq:magn} and \ref{eq:phase}. As an example, for a harmonic excitation of frequency $\omega$ = 10 rad/s, the value of $\alpha$ is calculated as 1.9903 - 0.0151i. Figure \ref{SS}a depicts the steady state response of both the integer order model and its equivalent fractional model, each obtained from their associated transfer functions, for $\omega$ = 10 rad/s. Both models return an equivalent steady state response, as expected due to the exact correspondence of the analytical transfer functions. Figure \ref{SS}b depicts the magnitude and phase of the I-SDOF and F-SDOF models across a range of frequencies. Both sets of magnitude and phase show that the transfer function of the F-SDOF model matches exactly the transfer function of the corresponding I-SDOF. It should be pointed out that the phase of the fractional models is actually the negative value of the phase of the integer models as evidenced by the previously derived equations of the steady state response.

Figure \ref{SS}c shows a plot of the value of $\alpha$ as a function of $\omega$. The natural frequency of both the I-SDOF and the F-SDOF in this numerical example is $\sqrt[]{5}$. Figure \ref{SS}c shows that for excitation frequencies above the critical frequency, the real part of $\alpha$ tends towards 2 while the imaginary part of $\alpha$ tends toward 0. Thus, when the forcing frequency is well above the resonance frequency (i.e. the high-frequency asymptotic limit), the fractional model converges to an integer order model. On the contrary, for frequencies below the critical frequency, the model order is fractional and $1<Re(\alpha)<2$ as expected due to the presence of viscous damping in the initial integer order model. 

\subsection{Model order reduction from I-MDOF to F-SDOF} \label{mdofsdof}

Considering the general procedure discussed above, we can discuss the case of model order reduction from a multiple integer-order degree of freedom (I-MDOF) system to a fractional SDOF. This is an extreme form of reduction where the initial MDOF system is collapsed to a single DOF model. This study case is chosen to show the flexibility of the fractional model approach and the ability to yield highly accurate, or even exact, solutions over an extended frequency range.

\begin{figure}[!ht] 
  \begin{center}
    \centerline{\includegraphics[scale=0.65,angle=0]{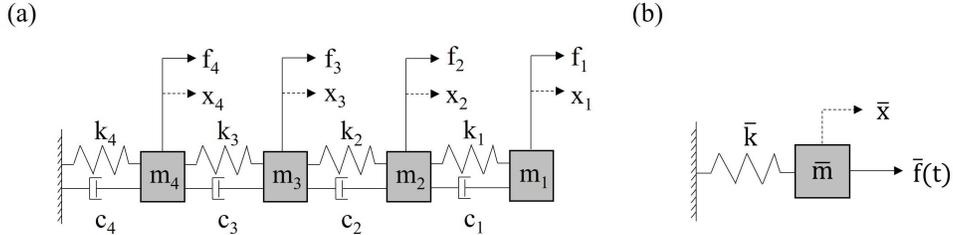}}
    \caption{Schematic of (a) the I-MDOF and of (b) the corresponding F-SDOF.}
    \label{D2}
  \end{center}
\end{figure}

Assume that the response of one of the degrees of freedom in the I-MDOF (Figure \ref{D2}a) model is the active degree of interest and that all the remaining DOFs are omitted. We seek an equivalent F-SDOF representation (Figure \ref{D2}b) of the system such that its response matches exactly or approximately the response of the corresponding active degree of freedom of the I-MDOF model. Recall that the transfer function of the fractional SDOF is given by Equation \ref{eq:TF}. In order to obtain the order $\alpha$ of the F-SDOF, we first convert the I-MDOF to a state-space form. From the state-space form, the transfer function of any of the nodes in the integer MDOF can be obtained. The process of obtaining the transfer function from state-space is well-established \cite{Franklin,Juang} and it is also briefly reviewed in \S \ref{App_B}. 

Let $H(s)$ be the transfer function for the active degree of freedom in the I-MDOF model, obtained from the state-space form. Assume that the corresponding fractional model has a mass of $\bar{m}$ and a stiffness of $\bar{k}$, as reflected in Equation \ref{eq:FRO} and Figure \ref{D2}b. By equating the transfer function of the fractional SDOF to the transfer function of the degree of interest in the I-MDOF, we obtain the fractional order $\alpha$:

\begin{eqnarray} \label{eq:alphaNDOF}
\alpha = \frac{\mbox{ln}(\frac{H_D}{\bar{m}H_N}-\frac{\bar{k}}{\bar{m}})}{\mbox{ln}(s)},
\end{eqnarray}

\noindent where $H$ is written as $H_N/H_D$. While we obtained $H$ by converting from state space to a transfer function, it should be noted that Equation \ref{eq:alphaNDOF} still holds true if the transfer function $H$ was obtained via another method. This latter characteristic will be the basis for the model synthesis approach discussed in \S \ref{Mod_syn}. Equation \ref{eq:alphaNDOF} can be made a function of frequency by substituting $s = i\omega$. Note that Equation \ref{eq:alphaNDOF} provides a fractional order $\alpha$ which guarantees an exact match of the transfer functions of the two systems. Therefore, as far as the individual functions $H_D$ and $H_N$ are exact, the response of the equivalent F-SDOF is an exact match of the initial I-MDOF.

Contrarily to the I-SDOF case, here, the equivalent mass and stiffness of the F-SDOF need to be determined. In fact, while in \S \ref{ISDOFtoFSDOF} the selection of the equivalent parameters $\bar{m}$ and $\bar{k}$ was straightforward given the existence of only one set of parameters, in the current configuration, multiple choices can be made. Among the possible approaches, we decided to set $\bar{m}$ so that the total mass of the I-MDOF matches that of the F-SDOF:
  
  \begin{eqnarray} \label{eq:mass}
  \bar{m} = \sum_{i=1}^{N} m_i,
  \end{eqnarray}
  
  \noindent where $m_i$ is the mass of the $i$-th degree of freedom of the I-MDOF. Since the I-MDOF model consists of springs in series, we set $\bar{k}$ equal to the equivalent stiffness of springs in series:
 
 \begin{eqnarray} \label{eq:stiffness}
 \bar{k} = \bigg(\sum_{i=1}^{N} \frac{1}{k_i}\bigg)^{-1}.
 \end{eqnarray}
 
\begin{figure}[!ht] 
  \begin{center}
    \centerline{\includegraphics[scale=0.7,angle=0]{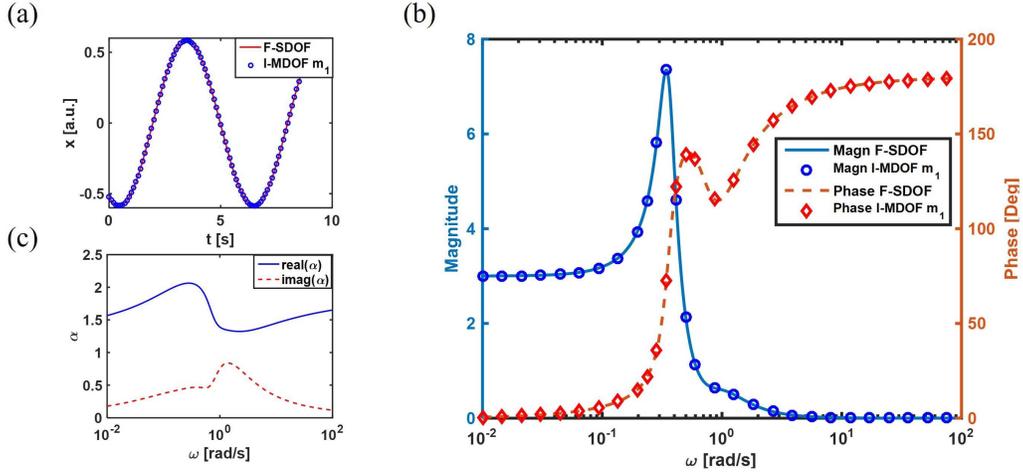}}
    \caption{a) Numerical simulations of the forced response showing the correspondence between the I-MDOF and the F-SDOF after the matching procedure is applied. Displacement time history at steady state for both the I-MDOF and the F-SDOF at a forcing frequency of 1 rad/s. b) Magnitude and phase of the transfer functions of I-MDOF and F-SDOF. c) Plot of the fractional order $\alpha(\omega)$ for the F-SDOF.}
    \label{NS}
  \end{center}
\end{figure}

As in the previous case, we numerically test the methodology. We consider an I-MDOF with $M=4$ and non-uniform coefficients. Specifically, we consider $m_1$ = $m_3$ = 1, $m_2$ = $m_4$ = 2, $k_1$ = $k_3$ = 1, $k_2$ = $k_4$ = 2, $c_1$ = $c_3$ =1, and $c_2$ = $c_4$ = 2. The natural frequencies are all within the range 0.36 rad/s to 2.22 rad/s. Let the dynamic response of the first mass ($m_1$ in Figure \ref{D2}a) be the active degree and therefore the one whose response is reduced to a fractional SDOF. Using Equations \ref{eq:mass} and \ref{eq:stiffness}, the parameters of the fractional model are taken as $\bar{m}$ = 6 and $\bar{k}$ = $\frac{1}{3}$. Assuming that we have obtained the transfer function $H(s)$ for the first mass in the integer MDOF model, Equation \ref{eq:alphaNDOF} can be applied to find $\alpha$ for a desired forcing frequency. The steady state response of the fractional model can be obtained using the magnitude and phase as given by Equations \ref{eq:magn} and \ref{eq:phase}. For a harmonic forcing function of frequency $\omega$ = 1 rad/s, the value of $\alpha$ is calculated as 1.3807 + 0.7731i. Figure \ref{NS}a shows the steady state response of both the first mass of the integer order model and its equivalent fractional model, each obtained from their associated transfer functions. As expected, the response of the equivalent fractional model is identical to the response of the first mass of the I-MDOF. Figure \ref{NS}b depicts the match between the magnitude and phase of the first mass of the I-MDOF and the magnitude and phase of the F-SDOF.

Recall that the transfer function $H$ is a function of the forcing frequency $\omega$; thus Equation \ref{eq:alphaNDOF} is a function of $\omega$.
 Figure \ref{NS}c shows a plot of the value of $\alpha$ as a function of $\omega$ for our example system. The trend of $\alpha(\omega)$ is highly dependent on the methodology used for determining the variables $\bar{m}$ and $\bar{k}$. Following our suggested methodology (Equations \ref{eq:mass} and \ref{eq:stiffness}) results in a trend of $\alpha(\omega)$ that closely follows that observed in \S \ref{mdofsdof}. In the asymptotic limit, the fractional system tends to a second order system ($Re(\alpha)=2$) with no phase modulation ($Im(\alpha)=0$); that is, no frequency dependence of the phase. Thus, in the asymptotic limit, the phase of the response becomes nearly constant and no longer changes with frequency. In the low frequency limit, we observe a fractional oscillator behavior ($1<Re(\alpha)<2$) with phase modulation ($Im(\alpha)\neq0$). Note that the original system is of order eight with resonances clustered in a narrow frequency range. For this reason, in the frequency range of the local resonances, the equivalent fractional system can exhibit order $Re(\alpha)>2$. As highlighted above, the fractional order is affected by the choice made to define the equivalent parameters. Results also confirm that non-monotonic changes in the fractional order should be expected in between resonances, suggesting that the evolution from one local resonance to the next (typically modeled in conventional dynamic theory as a SDOF second order oscillator) occurs via a dynamic behavior that is locally fractional.
 
 We conclude this section noting the two main advantages of the technique presented above: 1) the approach allows a remarkable reduction in order without any loss in the characteristic features of the original dynamics, and 2) the resulting F-SDOF can pave the way to the use of analytical solutions (e.g. Mittag-Leffler function \cite{Podlubny}) for the simulation of complex dynamic systems that typically allow only numerical approaches.

\section{Reduction to NDOF fractional models}

\begin{figure}[h] 
  \begin{center}
    \centerline{\includegraphics[scale=0.65,angle=0]{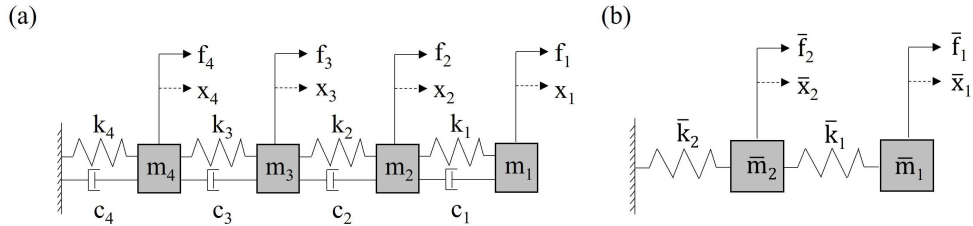}}
    \caption{Schematic of (a) the I-MDOF and of (b) the corresponding F-NDOF.}
    \label{D3}
  \end{center}
\end{figure}

In the previous section, we successfully reduced an I-MDOF to a F-SDOF while exactly retaining the dynamic response of the degree of interest in the initial integer order model. We extend this approach in order to make the reduction procedure more general and account for the reduction of an I-MDOF to a fractional multiple degree of freedom system (F-NDOF) with a number of DOFs $N<M$, such as the models depicted in Figure \ref{D3}.  This approach allows reducing the original I-MDOF model while retaining multiple active degrees of freedom.

\subsection{NDOF fractional Model}

We consider how to obtain the transfer functions for a F-NDOF from the equations of motion. The equations of motion of each degree in the F-NDOF are actually very similar to those of a I-MDOF, except for the fact that the inertia terms ($m \frac{d^2 x}{dt^2}$) are replaced by fractional terms ($\bar{m} \frac{d^\alpha x}{dt^\alpha}$). Recall that the selected fractional models in our approach do not include any explicit damping term. Consider a F-2DOF as in Figure \ref{D3}b. The equations of motion are given by:

\begin{eqnarray} \label{eq:FNDOF1}
\bar{m}_1\frac{d^\alpha \bar{x}_1}{dt^\alpha}+\bar{k}_1 \bar{x}_1 - \bar{k}_1 \bar{x}_2 = \bar{f}_1(t),
\end{eqnarray}
  
\begin{eqnarray} \label{eq:FNDOF2}
\bar{m}_2\frac{d^\alpha \bar{x}_2}{dt^\alpha}-\bar{k}_1 \bar{x}_1 + (\bar{k}_1+\bar{k}_2) \bar{x}_2 = \bar{f}_2(t).
\end{eqnarray} 

\noindent We could arrange these equations in a state-space form and then obtain the transfer functions in an analogous manner to what was done for the I-MDOF. However, rather than using a fractional state-space to transfer function method, we opt to calculate the Laplace transform of Equations \ref{eq:FNDOF1} and \ref{eq:FNDOF2} and then apply Cramer's rule. Note that this approach could also be applied to obtain the transfer functions of the I-MDOF. Applying the Laplace transform yields:

\begin{eqnarray} \label{eq:LTNDOF1}
\bar{m}_1 s^\alpha \bar{X}_1 + \bar{k}_1 \bar{X}_1 - \bar{k}_1 \bar{X}_2 = \bar{F}_1,
\end{eqnarray}
  
\begin{eqnarray} \label{eq:LTNDOF2}
\bar{m}_2 s^\alpha \bar{X}_2 -\bar{k}_1 \bar{X}_1 + (\bar{k}_1+\bar{k}_2) \bar{X}_2 = \bar{F}_2.
\end{eqnarray}

\noindent In matrix form:

\begin{eqnarray}
\begin{bmatrix} \bar{m}_1 s^ \alpha + \bar{k}_1 & -\bar{k}_1 \\ -\bar{k}_1 & \bar{m}_2 s^\alpha + \bar{k}_1 +\bar{k}_2 \end{bmatrix} \begin{bmatrix} \bar{X}_1 \\ \bar{X}_2 \end{bmatrix} = \begin{bmatrix} \bar{F}_1 \\ \bar{F}_2 \end{bmatrix}.
\end{eqnarray}

\noindent Using Cramer's rule, $\bar{X}_1(s)$ and $\bar{X}_2(s)$ are given by:

\begin{eqnarray} \label{eq:X1s}
\bar{X}_1(s) = \frac{D_1(s)}{D(s)},
\end{eqnarray} 

\begin{eqnarray} \label{eq:X2s}
\bar{X}_2(s) = \frac{D_2(s)}{D(s)},
\end{eqnarray} 

\noindent where:

\begin{eqnarray}
D_1(s) = 
\begin{vmatrix} \bar{F}_1(s) & -\bar{k}_1 \\ \bar{F}_2(s) & \bar{m}_2 s^\alpha + \bar{k}_1 +\bar{k}_2 \end{vmatrix},
\end{eqnarray}

\begin{eqnarray}
D_2(s) = 
\begin{vmatrix} \bar{m}_1 s^\alpha + \bar{k}_1 & \bar{F}_1(s) \\ -\bar{k}_1 & \bar{F}_2(s) \end{vmatrix},
\end{eqnarray}

\begin{eqnarray}
D(s) = 
\begin{vmatrix} \bar{m}_1 s^\alpha + \bar{k}_1 & -\bar{k}_1 \\ -\bar{k}_1 & \bar{m}_2 s^\alpha + \bar{k}_1 +\bar{k}_2 \end{vmatrix}.
\end{eqnarray}

The desired transfer functions can then be obtained using Equations \ref{eq:X1s} and \ref{eq:X2s}. As an example, if $\bar{F}_2(s) = 0$ and $\bar{F}_1(s) = \bar{F}(s)$ we can define $\hat{D}_1(s) = \frac{D_1(s)}{\bar{F}(s)}$ and $\hat{D}_2(s) = \frac{D_2(s)}{\bar{F}(s)}$. The transfer functions of the output displacements relative to the input force on the first node are:

\begin{eqnarray} \label{eq:F2DOF1}
G_1(s) = \frac{\bar{X}_1(s)}{\bar{F}(s)} = \frac{\hat{D}_1 (s)}{D(s)},
\end{eqnarray}

\begin{eqnarray} \label{eq:F2DOF2}
G_2(s) = \frac{\bar{X}_2(s)}{\bar{F}(s)} = \frac{\hat{D}_2 (s)}{D(s)}.
\end{eqnarray}

\noindent The substitution $s=i\omega$ is made to convert to the frequency domain. 

\subsection{Order reduction from I-MDOF to F-NDOF}

Leveraging the transfer function representation of the F-NDOF system provided in the previous paragraph, we can illustrate the reduction procedure from I-4DOF to F-2DOF. While, in general, both the I-MDOF's and the F-NDOF's degrees can have a forcing function acting on it, this example shall only include an external force on the first node in both the I-4DOF and in the F-2DOF. In the reduction process, let the 1st and 3rd masses in the I-4DOF be the active degrees. Therefore, the reduction technique shall yield steady state responses of the F-2DOF's degrees which are equivalent to the steady state responses of the active degrees in the I-4DOF. Assume we have obtained the transfer functions (that is, the transfer function between the displacement and the force) of active degrees in the integer order model by using state space to transfer function techniques. We refer to these transfer functions as $H_1$ and $H_3$, respectively. 

In order to match the responses of the two systems at the selected DOFs, we impose $H_1$ = $G_1$ and $H_3$ = $G_2$ where $G_1$ and $G_2$ are given by Equations \ref{eq:F2DOF1} and \ref{eq:F2DOF2}, respectively. It is tempting to arrange Equation \ref{eq:LTNDOF2} (with $\bar{F}_2$=0) into a form of $\frac{X_2}{X_1}$ and match it to $\frac{H_3}{H_1}$ and solve for $\alpha$. However, this will not lead to the desired matching of both responses of the active DOFs. To understand why, consider the mass and stiffness parameters in the F-2DOF. In the above procedure, it was assumed that the values of $\bar{m}_1$, $\bar{m}_2$, $\bar{k}_1$, and $\bar{k}_2$ were determined by lumping masses and stiffnesses from the integer model. If we set $G_1$ equal to $H_1$ and $G_2$ equal to $H_3$, we obtain a system of two equations in one unknown ($\alpha$), resulting in an over-determined system. Unless $\bar{m}_1$, $\bar{m}_2$, $\bar{k}_1$, and $\bar{k}_2$ are properly chosen, the resulting fractional order $\alpha$ will not be able to guarantee the exact matching of the reduced DOFs. Although different procedures to identify the parameters can be selected, we propose the following method. Lump the masses accordingly so that the total mass of the fractional model is equal to the total mass of the integer model. In our example, we set $\bar{m}_1$ = $m_1 + m_2$ and $\bar{m}_2$ = $m_3+m_4$. Next, set a relationship among the stiffness values of the F-NDOF. Recall Equation \ref{eq:stiffness}, in which we obtained a parameter $\bar{k}$ by finding the equivalent stiffness of springs in series. We define $\bar{k}_1$ = $\beta \bar{k}$ and $\bar{k}_2$ = $(1-\beta)\bar{k}$, where $\beta$ is an unknown coupling stiffness parameter. Note that $\bar{k}_1 + \bar{k}_2$ = $\bar{k}$. Equating the transfer functions of the fractional and integer order models, we obtain a set of two nonlinear equations in two unknowns ($\alpha$ and $\beta$). To solve for the complex order $\alpha$ and the stiffness coupling parameter $\beta$ (which will also be a complex quantity), a nonlinear numerical solver can be used. 

\begin{figure}[!ht] 
  \begin{center}
    \centerline{\includegraphics[scale=0.7,angle=0]{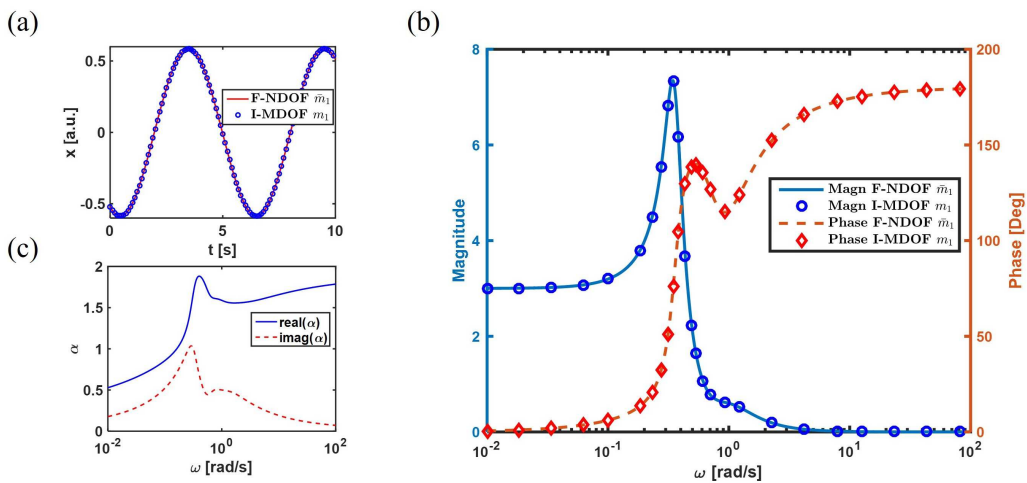}}
    \caption{a) Numerical simulations of the forced response showing the correspondence in response of $x_1$ in the I-MDOF model and $\bar{x}_1$ in the F-NDOF model. Displacement time history at steady state for both the I-MDOF and the F-NDOF at a forcing frequency of 1 rad/s. b) Magnitude and phase of the transfer functions of $m_1$ in the I-MDOF model and $\bar{m}_1$ in the F-NDOF model. c) Plot of the fractional order $\alpha(\omega)$ for the F-NDOF.}
    \label{NN}
  \end{center}
\end{figure}

\begin{figure}[!ht] 
  \begin{center}
    \centerline{\includegraphics[scale=0.7,angle=0]{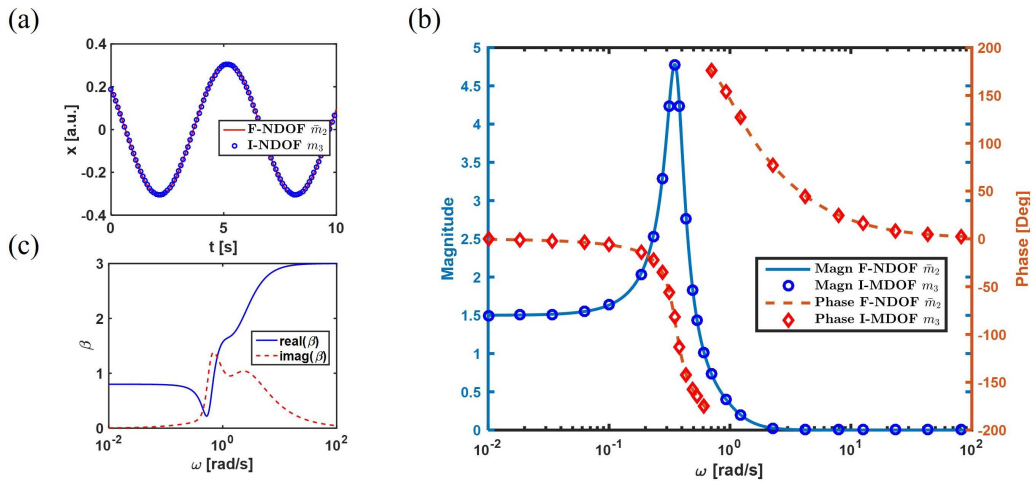}}
    \caption{a) Numerical simulations of the forced response showing the correspondence in response of $x_3$ in the I-MDOF model and $\bar{x}_2$ in the F-NDOF model. Displacement time history at steady state for both the I-MDOF and the F-NDOF at a forcing frequency of 1 rad/s. b) Magnitude and phase of the transfer functions of $m_3$ in the I-MDOF model and $\bar{m}_2$ in the F-NDOF model. c) Plot of the variable $\beta(\omega)$ for the F-NDOF.}
    \label{NNb}
  \end{center}
\end{figure}

As previously done, the procedure is further illustrated using a numerical example. Assume $m_1$ = $m_3$ = 1, $m_2$ = $m_4$ = 2, $k_1$ = $k_3$ = 1, $k_2$ = $k_4$ = 2, $c_1$ = $c_3$ = 1, and $c_2$ = $c_4$ = 2. The equivalent mass coefficients in the F-NDOF are $\bar{m}_1$ = $\bar{m}_2$ = 3. The magnitudes and phases over a frequency range are given in Figures \ref{NN}b and \ref{NNb}b. Figure \ref{NN}b depicts the magnitude and phase of the first degree (mass) in the I-4DOF and the first degree of the F-2DOF while Figure \ref{NNb}b show the magnitude and phase of the third degree in the I-4DOF and the second degree of the F-2NDOF. Results illustrate that the dynamic response of the F-NDOF was able to exactly match the response of the active degrees of the I-MDOF. 

Furthermore, we verify that the steady state responses of the fractional and integer models are equivalent for a specific forcing frequency. For a harmonic forcing function of frequency $\omega$ = 1 rad/s, we find $\alpha$ = 1.5834 + 0.4983i and $\beta$ = 1.5896 + 1.0440i. The numerical steady state results are plotted in Figures \ref{NN}a and \ref{NNb}a. Figure \ref{NN}a shows the steady state response of the first mass $x_1(t)$ in the I-4DOF and the first mass $\bar{x}_1(t)$ of the F-2DOF. Figure \ref{NNb}a depicts the steady state response of the third mass $x_3(t)$ in the I-4DOF and the second mass $\bar{x}_2(t)$ of the F-2DOF. The perfect overlap between the two plots shows that the dynamic response of the F-NDOF exactly match the response of the active degrees of the I-MDOF.

Finally, Figure \ref{NN}c shows a plot of the value of $\alpha$ as a function of $\omega$ and Figure \ref{NNb}c plots the value of $\beta$ as a function of $\omega$. Once again, we observe the familiar trend in Figure \ref{NN}c where the $Re(\alpha)\rightarrow2$ while $Im(\alpha)\rightarrow0$ in the asymptotic limit. A similar trend is also exhibited by the coupling parameter $\beta$. At frequencies well past the resonance frequencies, $\beta$ becomes purely real. For ranges in which the coupling parameter is complex, $Re(\beta)$ is representative of the stiffness in the I-4DOF while $Im(\beta)$ (along with the complex order $\alpha$) contributed to the damping. For high frequencies, where both $\alpha$ and $\beta$ are integers, damping no longer plays a significant role in the steady state dynamics.

\section{System identification from numerical data} \label{Mod_syn}

The formulation developed in the previous sections assumed that the initial I-MDOF system was provided and that its fractional form was sought. However, there are many situations of practical interest in the analysis of dynamical systems where measured experimental data is available and either the corresponding mathematical model or its coefficients are unknown; this class of problems is typically referred to as system identification. Over the years, a variety of system identification techniques have been proposed \cite{Juang,Isermann}. One of the most common approaches for vibration problems relies on matching second-order systems to individual resonances, therefore approximating the response of the system at resonance as a second order oscillator. This approach is also at the basis of the conventional half-power bandwidth method \cite{Rao} for damping estimation.

This approach has some important limitations. When the resonance frequencies of the MDOF system are too closely spaced, the local resonance is not well approximated by the single DOF oscillator. Also, as shown in the previous numerical results, in between resonances, the behavior of the system is typically fractional due to the coupling between two or more DOFs. This also explains why, when comparing numerical and experimental data, the largest discrepancies are often observed at frequencies off-resonance (regardless of the accuracy at resonance).

The category of system identification methods requires the selection of a dynamic model that is matched to the experimental data by properly tuning the model parameters. In this approach, the structure of the mathematical model (typically based on differential operators) is selected \textit{a priori} without any detailed insight in the true physical nature of the system. We have shown above that different operating regimes would require different models to achieve an accurate representation. Fractional models offer a much more general approach to the formulation of the equations of motion because they are capable of capturing in a single mathematical model a variety of physical mechanisms that would otherwise require multiple integer order models. The well-known change in the dynamic behavior of a system when transitioning from the low to the high frequency regime is a classical example of this phenomenon. In addition, it should be considered that typically the most appropriate dynamical model to describe a complex system is not known \textit{a priori}; therefore, the use of a fractional model would allow a general approach in which the system identification process is allowed to converge to the most appropriate form of the governing equations. Therefore, in some respect this approach would result in a model identification method. 

The above considerations suggest that fractional order models might provide a powerful methodology for the dynamic characterization of complex systems from measured data. Here below, we illustrate the approach first for an equivalent F-SDOF and then for a F-NDOF system. Note that, in the following we will refer to \textit{measured} data as the data used to synthesize the dynamical models. In practice, we generated this reference data numerically.

\subsection{Fractional SDOF}
According to the general approach followed above, we first consider how to determine the value of $\alpha$ across a desired frequency range in a fractional SDOF. Assume that numerical (or experimental) data describing the dynamics of a system has been obtained at a single location or DOF. Source data will be presented in the form of Bode plots. Let the magnitude of the measured Bode plot at the desired frequency be $M$ and the phase be $\psi$. Furthermore, assume that the total mass and stiffness of the system have been measured. Note that an estimation of the system's stiffness could be obtained from the amplitude of the transfer function in the low frequency range (i.e. from the quasi-static limit). The mass and stiffness of the F-SDOF system are $\bar{m}$ and $\bar{k}$, respectively. Using these quantities in Equation \ref{eq:magn} and \ref{eq:phase} we can solve for $\tau$ and $\xi$ at a specific frequency. Rearranging these equations yields: 

\begin{eqnarray} \label{eq:exp1}
\tau^2 + \xi^2 = \frac{1}{(M\bar{m})^2},
\end{eqnarray}

\begin{eqnarray} \label{eq:exp2}
\xi = \tau \mbox{tan}(-\psi).
\end{eqnarray}

\noindent Substituting Equation \ref{eq:exp2} into Equation \ref{eq:exp1} results in:

\begin{eqnarray} \label{eq:exp3}
\tau^2(1+\mbox{tan}^2(-\psi)) = \frac{1}{(M\bar{m})^2}.
\end{eqnarray}

\noindent Using the identity 1+tan$^2$($\phi$) = sec$^2$($\phi$), Equation \ref{eq:exp3} can be written as:

\begin{eqnarray} \label{eq:exp4}
\tau^2 = \frac{1}{(M\bar{m})^2 \mbox{sec}^2(-\psi)}.
\end{eqnarray}

\noindent Taking the positive root: 

\begin{eqnarray} \label{eq:exptau}
\tau = \frac{\mbox{cos}(-\psi)}{\bar{m}M}.
\end{eqnarray}

\noindent Furthermore:

\begin{eqnarray} \label{eq:expxi}
\xi = \frac{\mbox{sin}(-\psi)}{\bar{m}M}.
\end{eqnarray}

Once $\tau$ and $\xi$ are obtained for a certain frequency $\omega$, the nonlinear Equations \ref{eq:tau} and \ref{eq:xi} can be numerically solved to find the coefficients $a$ and $b$ in $\alpha$ = $a+ib$. This procedure can then be repeated over a range of frequencies to obtain the variable order $\alpha(\omega)$ for the corresponding F-SDOF model.

\begin{figure}[!ht] 
  \begin{center}
    \centerline{\includegraphics[scale=0.7,angle=0]{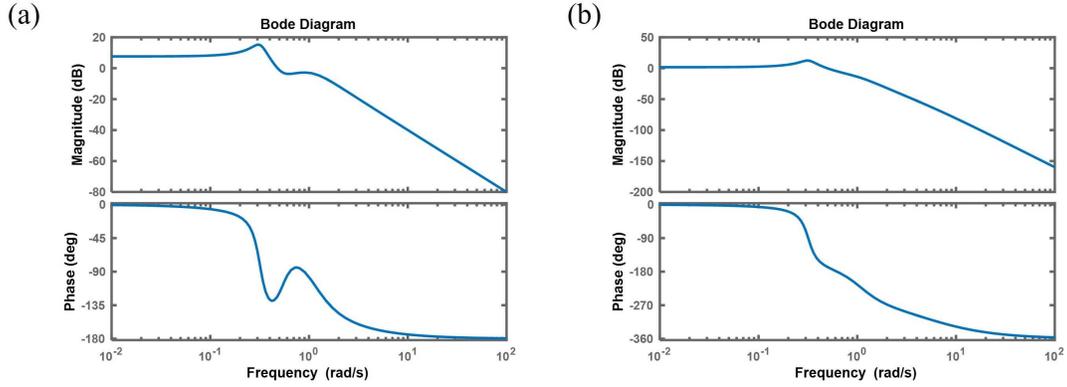}}
    \caption {a) The first set of measured magnitude and phase data. In an experiment setting, this could be understood as the data measured at a specific location.  b) The second set of measured magnitude and phase data. In an experiment setting, this could be understood as the data measured at another location different than that in part a). }
    \label{Bode}
  \end{center}
\end{figure}

\begin{figure}[!ht] 
  \begin{center}
    \centerline{\includegraphics[scale=0.7,angle=0]{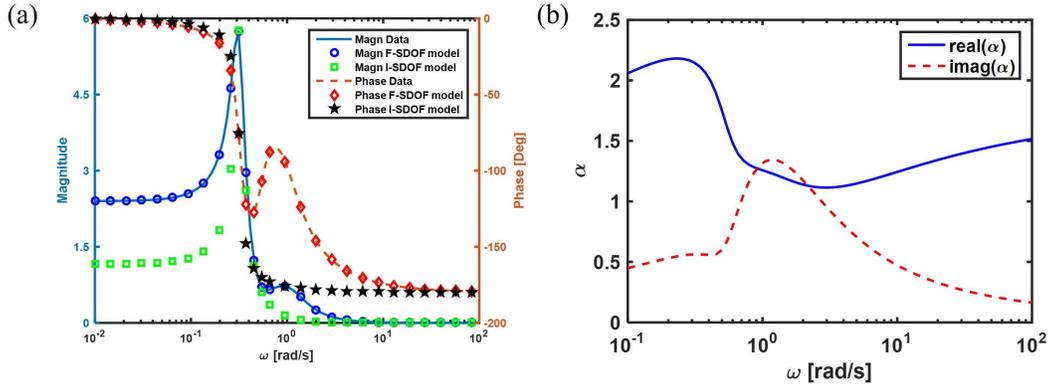}}
    \caption { a) The plot shows the magnitude and phase of measured data and of the transfer functions of the corresponding F-SDOF and I-SDOF models. b) Plot of the fractional order $\alpha(\omega)$ corresponding to the F-SDOF model synthesized from the measured data.}
    \label{EXPS}
  \end{center}
\end{figure}

To illustrate the procedure, assume that we have numerically (or experimentally) obtained the Bode plot given in Figure \ref{Bode}a. We wish to accurately represent this data using a F-SDOF model. Assume that the mass $\bar{m}$ and stiffness $\bar{k}$ of the system are known. In this example, we choose $\bar{m}$ = 12 and $\bar{k}$ = 0.4167. Equations \ref{eq:exptau} and \ref{eq:expxi} are solved to find $\tau$ and $\xi$ for 100 different frequencies between 0.01 and 100 rad/s. For each frequency, the value of $\alpha(\omega)$ can be found using Equations \ref{eq:tau} and \ref{eq:xi}. Figure \ref{EXPS}b shows the plot of $\alpha(\omega)$ for the selected example. To verify that the magnitude and phase of the F-SDOF match that of the numerical measured data, Equations \ref{eq:magn} and \ref{eq:phase} are used to plot the magnitude and phase of the obtained F-SDOF. Figure \ref{EXPS}a shows that the magnitude and phase of the F-SDOF match exactly the initial data. In order to assess the effectiveness of the fractional model approach, Figure \ref{EXPS}a also reports the magnitude and phase of an I-SDOF model created from the same input data. The parameters (mass, stiffness, damping) of the I-SDOF model were obtained by matching the transfer function of a second-order system at a selected resonance frequency ($\omega$ = 0.3 rad/s). Note that this approach follows the traditional method used to extract physical parameters (e.g. damping) from experimental data. In fact, for sufficiently spaced resonances, the peaks in the magnitude of the transfer function are fit locally by single DOF systems. While the obtained I-SDOF provides a good match of the data near the resonance, larger discrepancies are observable off-resonance. 

\subsection{Fractional NDOF}

Let us consider how to determine the value of $\alpha$ at a given frequency in a fractional NDOF. Recall in the fractional SDOF, we used data from one Bode plot (magnitude and phase) to find $\alpha$. This concept can be extended further. Consider that we have $N$ measured Bode plots obtained at different locations or degrees of interest and we want to represent the data by a fractional NDOF model. Furthermore, assume that the total mass and stiffness of the system have been measured. The masses in the fractional NDOF are determined by appropriately discretizing the total mass of the system into $N$ blocks. The stiffnesses of the F-NDOF are not yet known. However, the model will use \textit{stiffness coupling parameters} to establish relationships between the stiffnesses in the F-NDOF. There will be a total of $N$ unknowns; that is, the fractional order $\alpha$ and $N-1$ \textit{stiffness coupling} parameters.
If the analysis is conducted at a specific frequency, there are $2N$ known quantities from the experimental data; that is, the values of the magnitude and phase of each measured Bode plot. Considering the real and imaginary terms of $\alpha$ and the $N-1$ stiffness coupling parameters separately, there are $2N$ unknowns. Furthermore, there are $2N$ equations which are obtained by equating both the measured magnitude and phase of the transfer functions $G_N$ from the fractional model to the magnitude and phase of the available Bode plots. The set of $2N$ nonlinear equations can be solved numerically.

\begin{figure}[!ht] 
  \begin{center}
    \centerline{\includegraphics[scale=0.7,angle=0]{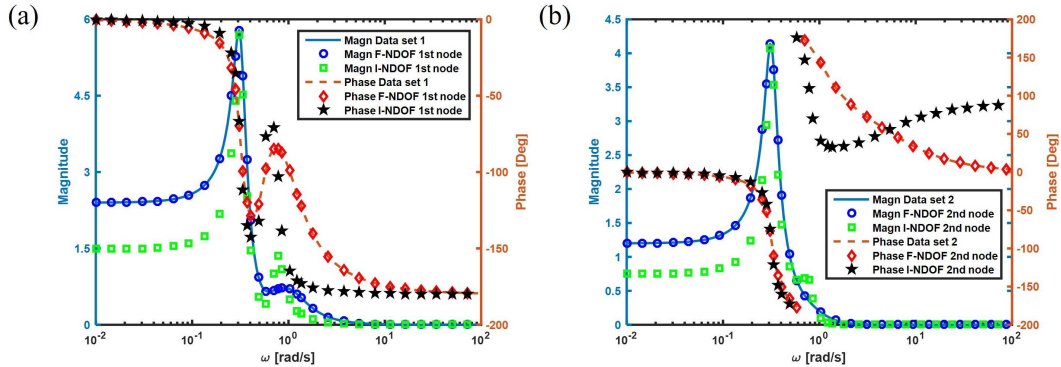}}
    \caption { Plot of magnitude and phase of the measured data in Bode form and of the transfer functions of the first masses in the corresponding F-2DOF and I-2DOF models. a) and b) show the results for the two different data sets.}
    \label{EXPNb}
  \end{center}
\end{figure}

\begin{figure}[!ht] 
  \begin{center}
    \centerline{\includegraphics[scale=0.7,angle=0]{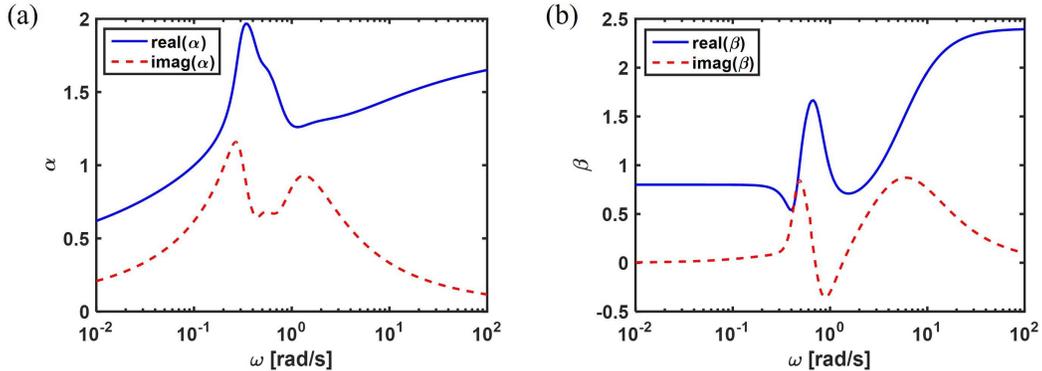}}
    \caption  { a) Plot of the fractional order $\alpha(\omega)$ for the corresponding F-2DOF model. b) Plot of the variable $\beta(\omega)$ for the corresponding F-2DOF model. }
    \label{EXPNa}
  \end{center}
\end{figure}

The procedure is illustrated on a F-2DOF model using two measured Bode plots (see Figures \ref{Bode}a and \ref{Bode}b) obtained at different locations. Assume that the total mass and stiffness are set to $\bar{m}$ = 12 and $\bar{k}$ = 0.4167 and that the corresponding parameters of the F-2DOF are $\bar{m_1}$ = $\bar{m_2}$ = 6. The unknown coupling parameter $\beta$ is set up such that $\bar{k}_1$ = $\beta \bar{k}$ and $\bar{k}_2$ = $(1-\beta)\bar{k}$. The analysis is performed for 100 different frequencies between the values of 0.01 and 100 rad/s. The results are presented in terms of the Bode plots of the transfer functions (Figure \ref{EXPNb}) which clearly show that the F-2DOf model very closely captures  the original dynamic behavior. Following an approach equivalent to the previous paragraph, we performed a system identification using also an integer order MDOF model. Figure \ref{EXPNb} shows the magnitude and phase of an I-2DOF model created from the same input data. The model parameters were obtained by matching the transfer functions of the second-order system to the magnitude and phase of the two sets of Bode plots at a selected resonance frequency ($\omega$ = 0.3 rad/s). As in the previous case, we observe that at resonance the response of the integer order system is able to capture the overall dynamics of the measured data, although without providing an exact fit. For frequencies off-resonance the discrepancies become more evident and the model fails in capturing the effective trend of the data.
The order $\alpha(\omega)$ of the corresponding F-2DOF model is plotted in Figure \ref{EXPNa}a while $\beta(\omega)$ is given in Figure \ref{EXPNa}b.

\section{Conclusion}

This study investigated the use of fractional order differential models to simulate the dynamic response of non-homogeneous discrete systems and to achieve efficient and accurate model order reduction. A limitation of many existing reduction order techniques, especially when applied to complex non-homogeneous systems, lies in their ability to accurately capture the local and global broadband dynamic response. We showed that the use of fractional differential equations having complex and variable order $\alpha(\omega)$ provides a viable and powerful alternative to conventional integer order models often guaranteeing, at least for the range of problems discussed in the present study, an exact match of the dynamic response. The frequency-dependent properties of the fractional operator enable frequency-dependent modulation of the phase and amplitude which are at the basis of the broad spectrum of problems that can be addressed with this type of models. Analytical and numerical results showed that the reduced fractional models can accurately (often times exactly) represent the dynamics of the active degrees of the initial system across a wide frequency spectrum. This is a key advantage of fractional models over traditional reduction techniques which typically yield narrowband performance and can guarantee only local accuracy. The proposed fractional approach was also tested in the frame of a system identification application. We argue that the use of fractional models does not only largely increase the accuracy of the parameter identification but it effectively results in a model identification approach. Numerical results showed that the fractional models were able to very accurately reflect the dynamics of the measured data, even at frequencies far off-resonance. This latter attribute is typically not achievable with conventional system identification techniques. 

In conclusion, this study laid the groundwork for model order reduction and system identification techniques based on complex, variable order fractional models. The methodology was shown to accurately represent the dynamics of non-homogeneous systems over a wide frequency spectrum.

\section{Appendix A: Basic Derivatives of Fractional Calculus}\label{App_A}

Several different definitions of a fractional derivative are available in the literature. Here below, we provide the two definitions that are related to this work: the Riemann-Liouville and the Caputo derivatives.

The of the Riemann-Liouville fractional derivative of order $\alpha$, for $t$ $\in$ [$a,b$] is: 

\begin{eqnarray}
D_{RL}^{\alpha} f(t) = \frac{1}{\Gamma(n-\alpha)}\frac{d^n}{dt^n}\int_{a}^{t} f(\tau)(t-\tau)^{n-\alpha -1} d\tau,
\end{eqnarray}

\noindent where $\Gamma(t)$ is the Gamma function and $n$=$\ceil{\alpha}$. On the other hand, the Caputo fractional definition is: 

\begin{eqnarray}
D_{C}^{\alpha} f(t) = \frac{1}{\Gamma(n-\alpha)}\int_{a}^{t} \frac{d^n f(\tau)}{d\tau^n} (t-\tau)^{n-\alpha -1} d\tau.
\end{eqnarray}

\noindent The Laplace Transform of Riemann-Liouville fractional derivative is:

\begin{eqnarray} \label{LTRL}
\int_{0}^{\infty} e^{-st} D_{RL}^\alpha f(t) dt = s^{\alpha}F(s) - \sum_{k=0}^{n-1} s^k [D_{RL}^{\alpha -k-1} f(t)]_{t=0}, \qquad n-1 < \alpha \leq n,
\end{eqnarray}

\noindent while the Laplace Transform of Caputo fractional derivative is:

\begin{eqnarray} \label{LTC}
\int_{0}^{\infty} e^{-st} D_C^\alpha f(t) dt = s^{\alpha}F(s) - \sum_{k=0}^{n-1} s^{\alpha -k-1} f^{(k)}(0), \qquad n-1 < \alpha \leq n.
\end{eqnarray}

From Equations \ref{LTRL} and \ref{LTC}, it is evident that the Laplace transform of the Caputo derivative uses the same initial values that a typical integer order problem does (first derivative, second derivative, etc). The initial values of the Riemann-Liouville definition are actually non-integer order derivative values of the function at $t=0$. The physical meaning of the necessary initial conditions using the Riemann-Liouville definition is an open question. On the other hand, the Caputo derivative lends itself to initial values which have a well-defined physical interpretation (initial position, velocity, acceleration, etc). Therefore, we shall use the Captuo fractional derivative in our analysis. 

\section{Appendix B: Brief Review of State Space to Transfer Function}\label{App_B}

A discrete MDOF system is represented by a set of $M$ second-order ordinary differential equations (ODE). To simplify their numerical solution, it is common to convert the original system into a system of $2M$ first-order ODEs. In dynamics, this approach is often referred to as the \textit{state-space} form. This can be represented by the vector equation:

\begin{eqnarray}
\frac{d\textbf{x}}{dt}=\textbf{Ax}+\textbf{B}u,
\end{eqnarray}

\noindent where the input is $u$ and the output is: 

\begin{eqnarray}
y = \textbf{Cx}+Du.
\end{eqnarray}

The matrix $\textbf{A}$ is the system matrix, $\textbf{B}$ the input vector, $\textbf{C}$ is the output vector, and D is a scalar called the direct transmission term. To illustrate, consider a 2-DOF mass-spring-damper system where there is forcing on the first block and the interested output is the displacement of the first mass. If $x_1$ and $x_2$ are the displacements of the two blocks, then the state-space is:

\begin{eqnarray}
\frac{d}{dt}\begin{bmatrix} x_1 \\ \dot{x}_1 \\x_2 \\ \dot{x}_2 \end{bmatrix} = 
\begin{bmatrix} 0 & 1 & 0 & 0  \\  -\frac{k_1}{m_1} & -\frac{c_1}{m_1} & \frac{k_1}{m_1} & \frac{c_1}{m_1} \\0 & 0 & 0 & 1  \\ \frac{k_1}{m_2} & \frac{c_1}{m_2} & -\frac{k_1+k_2}{m_2} & -\frac{c_1+c_2}{m_2} \end{bmatrix}\begin{bmatrix} x_1 \\ \dot{x}_1 \\x_2 \\ \dot{x}_2 \end{bmatrix} + \begin{bmatrix} 0 \\ \frac{1}{m_1} \\0 \\ 0 \end{bmatrix}f(t).
\end{eqnarray}

\noindent Thus:

\begin{eqnarray}
\textbf{A} = 
\begin{bmatrix} 0 & 1 & 0 & 0  \\  -\frac{k_1}{m_1} & -\frac{c_1}{m_1} & \frac{k_1}{m_1} & \frac{c_1}{m_1} \\0 & 0 & 0 & 1  \\ \frac{k_1}{m_2} & \frac{c_1}{m_2} & -\frac{k_1+k_2}{m_2} & -\frac{c_1+c_2}{m_2} \end{bmatrix},
\end{eqnarray}

\begin{eqnarray}
\textbf{B}= 
 \begin{bmatrix} 0 \\ \frac{1}{m_1} \\0 \\ 0 \end{bmatrix},
\end{eqnarray}

\begin{eqnarray}
\textbf{C}= 
 \begin{bmatrix} 1 & 0 & 0 & 0 \end{bmatrix},
\end{eqnarray}

\noindent and $D$ = 0 (as is usually the case). The transfer function of the desired degree is $H(s)$ [in this case, $H(s) = X_1(s)/F(s)$ since the first entry in $\textbf{C}$ is 1] and is given by: 

\begin{eqnarray}
 H(s) = C(sI-A)^{-1}B + D, 
\end{eqnarray}

\noindent where $I$ is the identity matrix of the appropriate size.


\end{document}